\documentclass[11pt]{amsart}

\usepackage{amstext,amssymb,amsmath,amsbsy}

\usepackage{amsmath,amssymb,latexsym,dsfont}
\usepackage[left=2.5cm,right=2.5cm,top=2.2cm,bottom=2.2cm]{geometry}

\usepackage{amsmath,amsfonts,amssymb,epsfig, textcomp}
\usepackage{graphicx,tikz}
\usepackage{hyperref, cleveref} 
\usepackage{cases,listings}
\usepackage{color} 
\usepackage{mathabx}
\usepackage{esint}

\newtheorem{theorem}{Theorem}
\newtheorem{lemma}{Lemma}

\newtheorem{proposition}{Proposition}
\newtheorem{corollary}{Corollary}
\newtheorem{remark}{Remark}

\newtheorem{definition}{Definition}
\theoremstyle{definition}
\numberwithin{equation}{section}
\numberwithin{theorem}{section}
\numberwithin{lemma}{section}
\numberwithin{proposition}{section}
\numberwithin{remark}{section}
\numberwithin{example}{section}
\numberwithin{corollary}{section}
\numberwithin{exercise}{section}
\numberwithin{definition}{section}

\newcommand{\dsp}{\displaystyle}

\newcommand{\sgn}{\mathrm{sgn}}

\newcommand{\supp}{\operatorname{supp}}

\newcommand{\eps}{\varepsilon}

\newcommand{\mC}{\mathbb{C}}
\newcommand{\mN}{\mathbb{N}}

\newcommand{\mS}{\mathbb{S}}
\newcommand{\mR}{\mathbb{R}}

\newcommand{\mZ}{\mathbb{Z}}

\newcommand{\CS}{\mathcal S}
\newcommand{\mc}{\mathrm{c}}

\newcommand{\hu}{\hat u}

\newcommand{\cF}{{\mathcal F}}

\newcommand{\cL}{{\mathcal L}}

\newcommand{\cS}{{\mathcal S}}

\newcommand{\tr}{^\text{tr}}
\newcommand{\R}{\mathbb{R}}

\newcommand{\be}{\begin{equation}}
\newcommand{\ee}{\end{equation}}

\newcommand{\ws}{\widetilde \sigma}

\newcommand{\hf}{\hat f}
\newcommand{\hg}{\hat g}
\newcommand{\ta}{\widetilde a}

\title[A multiplier theorem and its applications]{Multilinear multiplier theorems and their applications
\\
to the Jacobian and the Hessian determinant}

\author[H.-M. Nguyen]{Hoai-Minh Nguyen}
\author[B.  Perthame]{Beno\^it Perthame}

\address[H.-M. Nguyen]{Sorbonne Universit\'e, Universits\'e Paris Cit\'e, CNRS, INRIA, \newline
\indent Laboratoire Jacques-Louis Lions, LJLL, F-75005 Paris, France
}
\email{hoai-minh.nguyen@sorbonne-universite.fr}

\address[B.  Perthame]{Sorbonne Universit\'e, Universits\'e Paris Cit\'e, CNRS, INRIA, \newline
\indent Laboratoire Jacques-Louis Lions, LJLL, F-75005 Paris, France
}
\email{benoit.perthame@sorbonne-universite.fr} 

\begin{document}

\begin{abstract} We establish several variants of the multilinear multiplier theorem of Coifman and Meyer. We also present examples that are not covered by existing theories. Our motivation comes from applications to the definition of the Jacobian and Hessian determinant  in the distributional sense.
\end{abstract}

\maketitle 



\section{introduction}

Let $d \ge 1$ and $m \ge 2$, and let $\dsp a: (\mR^d)^m =  \underbrace{\mR^d \times  \cdots \times \mR^d}_{\text{$m$-times}} \to \mC$ be a bounded function. Consider the multilinear operator $T: (\cS(\mR^{d}) \big)^m \to  \mC$ defined by 
\be \label{def-T}
T (f_1, \dots, f_m) (x) = \int_{(\mR^{d})^m} a(\xi_1, \dots, \xi_m) \hf_1 (\xi_1) \cdots \hf_m (\xi_m) e^{i  x \cdot (\xi_1 + \dots + \xi_m)} \, d \xi_1 \cdots \, d \xi_m.
\ee
Here and in what follows, for $x, y \in \mR^d$, we denote by $x \cdot y$ their Euclidean scalar product, and 
for an appropriate function $f: \mR^d \to \mC$,  $\hat f$ denotes  its Fourier transform, i.e., 
$$
\hf(\xi) = \int_{\mR^d} f(x) e^{i \xi \cdot x} \, dx \quad \mbox{ for } \xi \in \mR^d. 
$$
Occasionally, we also use the notation $\cF (f)$ to denotes $\hf$.

Assume that, for all multi-indices $\alpha$  
\be \label{cond-CM}
|\nabla^{\alpha} a (\xi_1, \dots, \xi_m)| \le \frac{C_\alpha}{(|\xi_1| + \cdots +  |\xi_m|)^{|\alpha|}}, \qquad  \alpha| \le N, 
\ee
for some sufficiently large $N$, Coifman and Meyer  \cite{CM78-A} (see also \cite{CM78-C}), in connection with the study of commutators of singular integrals, initiated by Calderon~\cite{Calderon65}, proved the following well-known result. 

\begin{theorem}[Coifman \& Meyer] \label{thm-CM} 
For $d \ge 1$, $m \ge 2$, assume \eqref{cond-CM} and let $1 < p_1, \dots, \, p_m < + \infty$ and $1 \le r < + \infty$ be such that 
\be \label{thm-CM-rp}
\frac{1}{r} = \frac{1}{p_1} + \cdots + \frac{1}{p_m}. 
\ee
 Then 
\be\label{thm-CM-cl}
\|T(f_1, \dots, f_m)\|_{L^r(\mR^d)} \le C \|f_1 \|_{L^{p_1}(\mR^d)} \dots \| f_m\|_{L^{p_m} (\mR^d)},  
\ee
for some positive constant $C$ independent of $f_1$, \dots, $f_m$. 
\end{theorem} 
  
Assertion \eqref{thm-CM-cl} also holds in the case $m=1$. This is the classic result on the multiplier due to Mihlin \cite{Mihlin56} and H\"ormander \cite{Hormander60} in which one can take 
for $N$ any integer greater than $d/2$. In fact, they established the result under the following weaker condition on $a$ (for $m=1$): 
\be
\sup_{R > 0} \|a _R \|_{H^{n}\big(B_{\mR^d}(2) \setminus B_{\mR^d}(1) \big)} < + \infty, 
\ee
where $B_{\mR^d}(r)$ denotes the ball centered at 0 of radius $r$ of $\mR^d$ and
$$
a_R (\xi) = a(R \xi) \quad \mbox{ for } \xi \in \mR^d, 
$$
see, e.g., \cite[Theorem 7.9.5]{Hormander1}. 
It was shown by Tomita~\cite{Tomita10} for $r > 1$ and Grafakos and Si~\cite{GS12} for $r=1$ (see also  \cite[Theorem 7.5.5]{Grafakos2}) that one can replace \eqref{cond-CM} by the following 
\be \label{cond-CM-G}
\sup_{R > 0} \|a _R \|_{H^{mn}(B_{(\mR^d)^m}(2) \setminus B_{(\mR^d)^m}(1))} < + \infty, 
\ee
where $a_R(\xi_1, \dots, \xi_m) = a(R \xi_1, \dots, R \xi_m)$. This is consistent with Mihlin \& H\"ormander's result in the case $m = 1$.

\Cref{thm-CM} has been developed and extended in many directions such as the boundedness of the bilinear Hilbert transform \cite{LT97,LT99}, the multilinear  fractional integrals \cite{KS99}, the  Kato \& Ponce type inequalities \cite{KP88,MPTT04,MS13-1,MS13-2,GO14}. More topics and references can be found in \cite{MS13-1,MS13-2,Grafakos1,Grafakos2} and the references therein.

In this paper, we establish other variants of \Cref{thm-CM} and present some of their applications. The starting point of our work is as follows. Let $a_1, \dots, a_m: \mR^d \to \mC$ be bounded functions such that, for $j=1,\dots, m$, $\xi \in \mR^d\backslash \{0\}$ and multi-indices $\alpha \in \mN^d$
\be \label{cond-aj}
|\nabla^{\alpha} a_j (\xi)| \le \frac{C_\alpha}{|\xi|^{|\alpha|}} \quad  \mbox{ for } |\alpha| \le N, 
\ee  
for some sufficiently large $N$.  In particular, $a_j$ verifies Mihlin \& H\"ormander's condition. Define $a: (\mR^d)^m \to \mC$ by 
\be \label{def-a-motivation}
a(\xi_1, \dots, \xi_m) = a_1(\xi_1) \cdots a_m(\xi_m) \quad \mbox{ for } \xi_1, \dots, \xi_m \in \mR^d. 
\ee
Then $T$ defined in \eqref{def-T} with $a$ given in \eqref{def-a-motivation} can be written under the form 
$$
T(f_1, \dots, f_m) (x) = \frac{1}{(2 \pi)^{dm}} g_1(x) \cdots g_m(x) \quad \mbox{where} \quad  \hg_j (\xi) = a_j(\xi) \hf_j(\xi) \mbox{ for } \xi \in \mR^d. 
$$
Assume \eqref{thm-CM-rp}. Since, by Mihlin \& H\"ormander's multiplier theorem,  
$$
\| g_j\|_{L^{p_j}(\mR^d)} \le C \| f_j\|_{L^{p_j}(\mR^d)}, 
$$
it follows from H\"older's inequality that 
$$
\| T(f_1, \dots, f_m) \|_{L^r(\mR^d)} \le C \| f_1\|_{L^{p_1}(\mR^d)}  \cdots \| f_m\|_{L^{p_m}(\mR^d)}. 
$$
Note that $a$ given by \eqref{def-a-motivation}  satisfies the condition 
\be \label{cond-a-N}
|\nabla^{\alpha_1} \dots \nabla^{\alpha_m} a (\xi_1, \dots, \xi_m)| \le  \frac{C_{\alpha_1, \dots, \alpha_m}}{|\xi_1|^{|\alpha_1|} \cdots |\xi_m|^{|\alpha_m|}} \quad \mbox{ for } \xi_1, \dots, \xi_m \in \mR^d, 
\ee
for all multi-indices $(\alpha_1, \dots, \alpha_m) \in (\mN^d)^m$ with $|\alpha_j| \le N$  but this $a$ does not satisfy \eqref{cond-CM-G} in general. 

\medskip 
In this paper, we first establish a result of type  \eqref{thm-CM-cl} for which \eqref{cond-CM} is replaced by \eqref{cond-a-N}. To this end, we first introduce the following notation. 

\begin{definition} For $d \ge 1$, $m \ge 2$ let $\sigma$ be a function defined from  $(\mR^d)^m$  into $\mC$. The funciton $\sigma$ is said to be poly-homogeneous if 
$$
\sigma(t_1 \xi_1, \dots, t_m \xi_m) = \sigma(\xi_1, \dots, \xi_m) \mbox{ for } t_j > 0 \mbox{ and } \xi_j \in \mR^d \setminus \{0 \} \mbox{ with } j =1, \dots, m. 
$$
\end{definition}

Here are classes of  poly-homogeneous  functions. Let $\sigma_m: (\mR^d)^m \to \mC$ be a multilinear map, i.e., 
$$
\sigma_m \mbox{ is linear with respect to each variable $\xi_j \in \mR^d$ for $1 \le j \le m$}.  
$$
Set, for $\beta > 0$, 
\be \label{example1}
\sigma (\xi_1, \dots, \xi_m) = \frac{\big|\sigma_m(\xi_1, \dots, \xi_m)\big|^{\beta}}{|\xi_1|^\beta \cdots |\xi_m|^\beta} \quad \mbox{ for } \xi_1, \dots, \xi_m \in \mR^d. 
\ee
or for  $k \in \mN$, set 
\be \label{example2}
\sigma (\xi_1, \dots, \xi_m) = \frac{\big(\sigma_m(\xi_1, \dots, \xi_m) \big)^{k}}{|\xi_1|^k \cdots |\xi_m|^k} \quad \mbox{ for } \xi_1, \dots, \xi_m \in \mR^d. 
\ee
Then $\sigma$ is also poly-homogeneous of degree 0. For example, let $m =2$ and $\beta > 0$, then 
$$
\sigma (\xi_1, \xi_2) = \frac{(\xi_1 \cdot \xi_2 )^\beta}{|\xi_1|^\beta |\xi_2|^\beta} \quad  \mbox{ for } \xi_1, \xi_2 \in \mR^d, 
$$
is poly-homogeneous of degree 0. In particular, when  $m = d$ and $k \ge 0$, then 
$$
\sigma (\xi_1, \dots, \xi_d) = \frac{\big(\det(\xi_1, \dots, \xi_d)\big)^k}{|\xi_1|^k \cdots |\xi_d|^k}  \quad  \mbox{ for } \xi_1, \dots, \xi_d \in \mR^d
$$
is poly-homogeneous of degree 0. 

\medskip 
Our result in this direction, which proof is given in \Cref{sect-thm-Multiplier},  is the following.   

\begin{theorem} \label{thm-Multiplier} For $d \ge 1$, $m \ge 2$, let $\ta, \sigma: (\mR^d)^m \to \mC$ be such that 
\begin{itemize}
\item[a)]  $\sigma$ is poly-homogeneous of degree $0$. 

\item[b)] $\sigma \in H^{mN}((\mS^{d-1})^m)$ for $N  = n + 2 d$ with some integer $n > d/2$ \footnote{$\mS^{d-1}$ denotes the unit sphere of $\mR^d$.}.

\item[c)] $\ta \in  H^{mN} ((\mR^d)^m)$ verifies the condition \eqref{cond-CM-G} where $a$ is replaced by $\ta$. 

\end{itemize}
Let $1 < p_1, \dots, \, p_m < + \infty$ and $1 \le r < + \infty$ be such that 
\be
\frac{1}{r} = \frac{1}{p_1} + \cdots + \frac{1}{p_m}. 
\ee
Then, with $T$ is defined by \eqref{def-T} with $a = \ta \sigma$, we have
\be
\| T (f_1, \dots, f_m) \|_{L^r(\mR^d)} \le C \|f_1 \|_{L^{p_1}(\mR^d)} \cdots \|f_m \|_{L^{p_m} (\mR^d)}, 
\ee
Here $C$ denotes a positive constant independent of $f_1, \dots, f_m$. 
\end{theorem}

Here is a direct consequences of \Cref{thm-Multiplier}. 

\begin{corollary} Let $d \ge 1$, $m \ge 2$. Given $k\in \mN$, define $\sigma$ by \eqref{example2}, and given $\beta > 0$, define $\sigma$ by \eqref{example1}. Assume that $\beta$ is an integer or $\beta \ge m N$ with $N  = n + 2d$ with some integer $n> d/2$, which is not necessary to be an integer. 
Let $1 < p_1, \dots, \, p_m < + \infty$ and $1 \le r < + \infty$ be such that 
$$
\frac{1}{r} = \frac{1}{p_1} + \cdots + \frac{1}{p_m}. 
$$
Then, with $T$ being defined by \eqref{def-T} for $a = \sigma$, it holds 
\be
\| T (f_1, \dots, f_m) \|_{L^r(\mR^d)} \le C \|f_1 \|_{L^{p_1}(\mR^d)} \cdots \|f_m \|_{L^{p_m} (\mR^d)}, 
\ee
for some positive constant $C$ independent of $f_1$, \dots,  $f_m$. 
\end{corollary}

\Cref{thm-Multiplier} in the case $\ta \equiv 1$ is a variant of \Cref{thm-CM} for which the condition \eqref{cond-CM} is replaced by \eqref{cond-a-N} under the additional assumption that $a$ is poly-homogeneous of degree $0$. It is worth noting that condition \eqref{cond-a-N} alone is not sufficient to guarantee the conclusion, as shown in \cite{GK01}.

\medskip 
We next discuss another related multilinear multiplier result with \Cref{thm-Multiplier}, which proof is given in \Cref{sect-thm-Multiplier-2}. This result is motivated from the study of distributional Jacobian and Hessian determinant, which are discussed below.

\begin{theorem} \label{thm-Multiplier-2} Let $d \ge 2$ and $k, m \in \mN$. Assume that $2 \le m \le d$ and let $\sigma_m: (\mR^d)^m \to \mC$ be an alternating multilinear map, i.e., $\sigma_m$ is a multilinear map such that $\sigma_m(\xi_1, \cdots, \xi_m) = 0$ if $\xi_{j_1} = \xi_{j_2}$ for some $1 \le j_1 \neq j_2 \le m$.  
Let $\sigma: (\mR^d)^m \to \mC$ be defined by 
\be \label{thm-Multiplier2-sigma}
\sigma (\xi_1, \dots,  \xi_m) = \big(\sigma_m(\xi_1, \dots, \xi_m) \big)^k  \quad \mbox{ for } \xi_1, \dots, \xi_m \in \mR^d. 
\ee
Let $1 < p_1, \dots, \, p_m < + \infty$ and $1 \le r < + \infty$ be such that 
$$
\frac{1}{r} = \frac{1}{p_1} + \cdots + \frac{1}{p_m}. 
$$
Set 
\be
s = \frac{k (m-1)}{m} \quad \mbox{ and } \quad r_* = \frac{r}{r-1}. 
\ee
Then, it holds, for  $\varphi \in C^\infty_c(\mR^d)$, 
\be
\Big| \langle T (f_1, \dots, f_d), \varphi  \rangle \Big| \le C \|f_1 \|_{\cL^{p_1}_s(\mR^d)} \dots \|f_m \|_{\cL^{p_m}_s(\mR^d)} \|\varphi \|_{W^{k, r_{*}}(\mR^d)} .     
\ee 
\end{theorem}

For $s \ge 0$ and $1 < p < + \infty$, we denote by $\cL^{p}_s$ the Bessel potential space
\be
\cL^{p}_s = \Big\{f \in L^p(\mR^d); g \in L^p(\mR^d) \quad \mbox{where} \quad  \hg(\xi) = (1 + |\xi|^2)^{1/2} \hf (\xi) \mbox{ for } \xi \in \mR^d \Big\},  
\ee
and the corresponding norm
$$
\| f\|_{\cL^{p}_s} = \| g\|_{L^p(\mR^d)} \quad  \mbox{where} \quad  \hg(\xi) = (1 + |\xi|^2)^{1/2} \hf (\xi) \mbox{ for } \xi \in \mR^d. 
$$
It is well-known that, see, e.g., \cite[Chapter V]{Stein-70}, 
\be \label{Bessel-p1}
\cL^{p}_s (\mR^d) = W^{s, p} (\mR^d) \mbox{ for $s \in \mN$}, 
\ee
\be\label{Bessel-p2}
\cL^{p}_s (\mR^d) = W^{s, p} (\mR^d) \quad \mbox{ for } p \ge 2, s \ge 0, 
\ee
and 
\be \label{Bessel-p3}
W^{s, p} (\mR^d) \subset \cL^p_s (\mR^d)\quad  \mbox{ for } 1 < p < 2, s \ge 0. 
\ee
Here,  for $s \ge 0$ and $1 < p < + \infty$, $W^{s, p}(\mR^d)$ denotes the fractional Sobolev spaces. It is also known that $\cL^p_s$ is a subset of the Besov space $B^{s}_{p, 2}$ for $1 < p < 2$, see, e.g., \cite[Chapter V]{Stein-70}. 

\medskip 
Here are two applications of \Cref{thm-Multiplier-2}. The first one is to the definition and the properties of the Jacobian in the distributional sense, due to Brezis and Nguyen \cite{BrNg11}. To connect their definition, we need the following lemma, whose proof is elementary and omitted. 
\begin{lemma} \label{lem-det-J} Let $d \ge 2$ and let $u \in C^\infty_{\mc} (\mR^d, \mR^d)$. We have 
\begin{multline}
(i)^{d} \cF (\det \nabla u) (\xi)
=\int_{(\mR^d)^{d-1}} \sigma(\xi - \xi_1, \xi_1 - \xi_2, \cdots,  \xi_{d-2} - \xi_{d-1}, \xi_{d-1}) \times  \\[6pt]
\times \hu (\xi - \xi_{1}) \hu(\xi_1) - \xi_2)  \cdots \hu(\xi_{d-2} - \xi_{d-1}) \hu(\xi_{d-1}) \, d \xi_1 \dots  \, d  \xi_{d-1}, 
\end{multline}
where 
\be
\sigma(\xi_1, \dots,   \xi_d)= \det (\xi_1, \dots,   \xi_d) \quad \mbox{ for } \xi_1, \dots, \xi_d \in \mR^d. 
\ee
\end{lemma}

Recall that $\cF (\det \nabla u) $ denotes the Fourier transform of $\det \nabla u$. 

\medskip 
Applying \Cref{thm-Multiplier-2} and using \Cref{lem-det-J} with $r=1$, we obtain the following result. 

\begin{proposition} \label{pro-Jacobian} 
For $d \ge 2$, set $s = 1 - 1/ d$. Let $1 < p_1, \dots, \, p_d < + \infty$ be such that 
\be
\frac{1}{p_1} + \cdots + \frac{1}{p_d} = 1. 
\ee
For $u, v \in C^\infty_{\mc} (\mR^d, \mR^d)$ and $\varphi \in C^\infty_{\mc}(\mR^d, \mR)$, we have
\begin{multline}
\left| \int_{\mR^d} \Big( \det (\nabla u) - \det (\nabla v) \Big) \varphi  \, dx \right| \\[6pt] \le C \left( \prod_{k=1}^d \| u_k\|_{\cL^{p_k}_{s} (\mR^d)}  + \prod_{k=1}^d \| v_k\|_{\cL^{p_k}_{s} (\mR^d)} \right) \sum_{j=1}^d \frac{\| u_j - v_j \|_{\cL^{p_j}_{s} (\mR^d)}}{\| u_j\|_{\cL^{p_j}_{s} (\mR^d)} + \| v_j\|_{\cL^{p_j}_{s} (\mR^d)}}  \| \nabla \varphi  \|_{L^{\infty}(\mR^d)},    
\end{multline}
for some positive constant $C$ depending only on $d$, $p_1$, \dots, $p_d$. 
In particular, we have
\be\label{pro-Jacobian-cl2}
\left| \int_{\mR^d} \det (\nabla u)  \varphi  \, dx \right| \\[6pt] \le C \prod_{k=1}^d \| u_k\|_{\cL^{p_k}_{s} (\mR^d)}   \| \nabla \varphi  \|_{L^{\infty}(\mR^d)}.     
\ee
\end{proposition}

As a consequence of \Cref{pro-Jacobian}, one can define and estimate the distributional Jacobian for  maps $u = (u_1, \cdots, u_d) \in \cL^{p_1}_s (\mR^d) \times \cL^{p_d}_s (\mR^d)$ with $1 <  p_j < + \infty$, $\sum_{j=1}^d \frac{1}{p_j} = 1$, and $s = 1- 1/d$. Moreover, \eqref{pro-Jacobian-cl2} holds.

\medskip 
The second application of \Cref{thm-Multiplier-2} is to the definition of the Hessian determinant  in the distributional sense. To this end, we establish the following result. 

\begin{lemma} \label{lem-det-H} Let $d \ge 2$. We have 
\begin{multline}
(-1)^{d} d! \cF (\det \nabla^2 u) (\xi) 
=\int_{(\mR^d)^{d-1}} \sigma(\xi - \xi_1, \xi_1 - \xi_2, \cdots,  \xi_{d-2} - \xi_{d-1}, \xi_{d-1}) \times  \\[6pt]
\times \hu (\xi - \xi_1) \hu(\xi_1 - \xi_2)  \cdots \hu(\xi_{d-2} - \xi_{d-1}) \hu(\xi_{d-1}) \, d \xi_1 \dots  \, d  \xi_{d-1} , 
\end{multline}
where 
\be
\sigma(\xi_1, \dots,   \xi_d)= \big(\det (\xi_1, \dots,   \xi_d) \big)^2 \quad  \mbox{ for } \xi_1, \dots, \xi_d \in \mR^d . 
\ee
\end{lemma}

Applying \Cref{thm-Multiplier-2} and using \Cref{lem-det-H} with $r=1$, we obtain the following result.

\begin{proposition} \label{pro-Hessian} 
Let $d \ge 3$. Denote $s = 2 - 2/ d$ and let $1 < p_1, \dots, \, p_d < + \infty$ be such that 
$$\frac{1}{p_1} + \cdots + \frac{1}{p_d} = 1. $$ 
For $u, v \in C^\infty_{\mc} (\mR^d, \mR)$ and $\varphi \in C^1_{\mc}(\mR^d, \mR)$, we have
\begin{multline}
\left| \int_{\mR^d} \Big( \det (\nabla^2 u) - \det (\nabla^2 v) \Big) \varphi  \, dx \right| \\[6pt] \le C \left( \prod_{k=1}^d \| u\|_{\cL^{p_k}_{s} (\mR^d)}  + \prod_{k=1}^d \| v\|_{\cL^{p_k}_{s} (\mR^d)} \right) \sum_{j=1}^d \frac{\| u - v \|_{\cL^{p_j}_{s} (\mR^d)}}{\| u\|_{\cL^{p_j}_{s} (\mR^d)} + \| v\|_{\cL^{p_j}_{s} (\mR^d)}}   \| \nabla^2 \varphi  \|_{L^{\infty}(\mR^d)},    \end{multline}
 for some positive constant $C$ depending only on $d$, $p_1$, \dots, $p_d$. 
In particular, 
\be\label{pro-Hessian-cl2}\left| \int_{\mR^d} \det (\nabla^2 u)  \varphi  \, dx \right| \\[6pt] \le C \prod_{k=1}^d \| u\|_{\cL^{p_k}_{s} (\mR^d)}   \| \nabla^2 \varphi  \|_{L^{\infty}(\mR^d)}.     \ee

\end{proposition}

\begin{remark}  \rm The result given in \Cref{pro-Hessian} also holds for $d=2$. Nevertheless, in this case the Hessian determinant has a simple formula and one can derive these properties easily, see, e.g., \Cref{rem-Hessian}. 
\end{remark}

It would be interesting to explore the connection of the results given in \Cref{thm-Multiplier-2} with the study of the Monge-Ampere equations, see, e.g., \cite{Brenier91,Gutierrez01,Figalli17,Le24}. 

\medskip 
We next briefly describe the idea of the proof of \Cref{thm-Multiplier,thm-Multiplier-2}. The proof of \Cref{thm-Multiplier} is inspired by the work of Coifman and Meyer \cite{CM78-A}. The idea is to localise each frequency variable $\xi_j$ in an annulus of the size $2^{k_j}$ for $1 \le j \le m$. This is nevertheless different with the proof of \Cref{thm-CM} where we first localise the vector frequency variable $(\xi_1, \cdots, \xi_m)$ in an annulus of the size $2^k$. Then, we introduce  the function $\sigma_{j_1, \dots, j_m} (\xi_1, \dots, \xi_m))$ in $\Big(B_{\mR^d} (2) \setminus B_{\mR^d}(1) \Big)^m$ defined by $\sigma(2^{j_1} \xi_1, \dots, 2^{j_m} \xi_m)$  which we decompose into an infinite sum which terms are of the form  $\varphi_1(\xi_1) \cdots \varphi_m(\xi_m)$. Applying \Cref{thm-Multiplier}, one can show that the desired conclusion holds for each term. The key observation is that one can obtain the corresponding estimate for the infinite sum of these terms thanks to the poly-homogeneity of degree $0$ of $\sigma$. The details are provided in \Cref{sect-thm-Multiplier}. 

The proof of \Cref{thm-Multiplier-2} uses \Cref{thm-Multiplier}. Nevertheless, the key argument is to localise the frequency of $(\xi_1, \dots, \xi_m)$ into a small cone  where there exists $1 \le j_0 \le m$ such that 
$$
| \xi_{j_0}| \le C |(\xi_1, \cdots, \xi_m)|. 
$$
The idea is then to write $\sigma$ given by \eqref{thm-Multiplier2-sigma} in the form of the multiplier given in \Cref{thm-Multiplier} and derive the corresponding result from it.  It is the place where the multilinear alternating property of $\sigma_m$ plays an important role. The details are given in \Cref{sect-thm-Multiplier-2}.

We next make some comments on the known works on distributional definitions of the Jacobian and Hessian determinant. It is wellknown, see  Morrey~\cite{Morrey}, Reshetnyak~\cite{Reshentnyak}, and Ball~\cite{Ball-convex}, that for $u\in C^1(\R^d; \R^d)$
\be\label{piola1}
\operatorname{det}(\nabla u) = \sum_{i=1}^d \frac{\partial}{\partial x_i} ( u_j C_{ij})  \quad  \mbox{ for } j=1,\dots, d,
\ee
where $C_{ij}$ denotes the  matrix of cofactors of the matrix $\nabla u)$. This directly follows from  the Piola idendity
\[
\sum_{i=1}^d \frac{\partial}{\partial x_i} C_{ij}  = 0 \quad  \mbox{ for } j=1,\dots, d.
\]
This mades it possible to define the distributional Jacobian determinant $\operatorname{Det}(\nabla u)$ for appropriate class of differentiability and integrability of $u$, which is larger than $W^{1, d}(\mR^d)$. Nevertheless, $u$ is required to have one derivative in the Sobolev sense in these works. The definition and the estimate of the distributional Jacobian in the fractional Sobolev scale was proposed and established by Brezis and Nguyen \cite{BrNg11} (see also \cite{BrNgJacobian} for a related context). They proved that one can define the distributional Jacobian for maps $u \in \big(W^{s, d}(\mR^d) \big)^d$ with $s = 1 - 1/d$. Moreover, it holds
$$
\left| \int_{\mR^d} \operatorname{det}(\nabla u) \psi  \right|
\le C \| u \|_{W^{s, d}(\mR^d)}^d \| \nabla \psi  \|_{L^\infty(\mR^d)}. 
$$
Their analysis is based on the integration by parts in $\mR^{d+1}_+$, one dimension thus is added. Their analysis also gives \eqref{pro-Jacobian-cl2} in which the space $\cL^{p_j}_s(\mR^d)$ is replaced by $W^{s, p_j}(\mR^d)$. As mentioned in \eqref{Bessel-p2} that $\cL^{p}_s(\mR^d) = W^{s, p}(\mR^d)$  when $p \ge 2$,  and in 
\eqref{Bessel-p3} that $W^{s, p}(\mR^d) \subset \cL^{p}_s(\mR^d)$ when $1 < p < 2$. We thus improve slightly their result using a completely different technique. 

Baer and Jerison \cite{BJ15} extend the analysis of Brezis and Nguyen to the distributional Hessian. They showed that one can define and estimate the Hessian determinant in the distributional sense. More precisely, they established, for $d \ge 2$, with $s = 2 - 2/d$, 
$$
\left| \int_{\mR^d} \operatorname{det}(\nabla^2 u) \psi \right|
\le C \| u \|_{W^{s, d}(\mR^d)}^d  \| \nabla^2 \psi \|_{L^{\infty}(\mR^d)}.  
$$
Their analysis can be extended to reach \eqref{pro-Hessian-cl2} in which the space $\cL^{p_j}_s(\mR^d)$ is replaced by $W^{s, p_j}(\mR^d)$. 
As mentioned in \eqref{Bessel-p2} that $\cL^{p}_s(\mR^d) = W^{s, p}(\mR^d)$ when $p \ge 2$,  and in 
\eqref{Bessel-p3} that $W^{s, p}(\mR^d) \subset \cL^{p}_s(\mR^d)$ when $1 < p < 2$. We thus improve lightly their result using a completely different technique. Notice that a formula analogous to \eqref{piola1} holds for the Hessian determinant. See Remark~\ref{Hessian_formula}. 

The paper is organized as follows. The proof of \Cref{thm-Multiplier,thm-Multiplier-2} are given in \Cref{sect-thm-Multiplier,sect-thm-Multiplier-2}, respectively. Finally, the proof of \Cref{lem-det-H} is given in \Cref{sect-lem-det-H}.

\section{Proof of \Cref{thm-Multiplier}} \label{sect-thm-Multiplier}

 For notational ease, we only present the proof in the case $m=2$. The proof in the general case follows similarly. We thus assume that $m=2$ later on in this proof. 
 
 Let $\psi \in C^\infty (\mR^d)$ be such that 
$$
\supp \psi \subset B_2 \setminus B_{1/2}, 
$$
and 
$$
\sum_{j \in \mZ} \psi (2^{j} \cdot) = 1 \mbox{ in } \mR^d. 
$$
We have 
$$
 \sigma(\xi_1, \xi_2) \hat f_1 (\xi_1)  \hat f_2 (\xi_2) = \sum_{j, k \in \mZ} \sigma(\xi_1, \xi_2)  \hat f_{1, j} (\xi_1) \hat f_{2, k} (\xi_2) \mbox{  for $\xi_1, \xi_2 \in \mR^d$,}
$$
where, for $j, k \in \mZ$,   
$$
\hat f_{1, j} (\xi_1)  = \psi (2^{-j} \xi_1) \hat  f_{1} (\xi_1)  
\quad \mbox{ and } \quad 
\hat f_{2, k} (\xi_2)  = \psi (2^{-k} \xi_2) \hat  f_{2} (\xi_2)  . 
$$

Let $\phi \in C^\infty(\mR^d)$ be such that 
$$
\supp \phi \subset B_{4} \setminus B_{1/4} \quad \mbox{ and } \quad \phi = 1 \mbox{ in } B_2 \setminus B_{1/2}. 
$$
Set, for $(\xi_1, \xi_2) \in (\mR^d)^2$, 
\be \label{thm-Multiplier-sigma-jk}
\widetilde \sigma_{j, k} (\xi_1, \xi_2) = \phi (\xi_1) \phi(\xi_2)  \sigma(2^j \xi_1, 2^k \xi_2).
\ee

We then have 
\begin{multline*}
\sigma(\xi_1, \xi_2) \hat f_1 (\xi_1)  \hat f_2 (\xi_2)  = 
\sum_{j, k \in \mZ} \sigma (\xi_1, \xi_2) \hat f_{1, j} (\xi_1) \hat f_{2, k} (\xi_2) \\[6pt]
= \sum_{j, k \in \mZ} \phi (2^{-j} \xi_1) \phi(2^{-k} \xi_2)  \sigma (\xi_1, \xi_2) \hat f_{1, j} (\xi_1) \hat f_{2, k} (\xi_2)\\[6pt]
= \sum_{j, k \in \mZ} \widetilde \sigma_{j, k} (2^{-j} \xi_1, 2^{-k} \xi_2) \hat f_{1, j} (\xi_1) \hat f_{2, k} (\xi_2)  \quad \mbox{ for } (\xi_1, \xi_2) \in (\mR^d)^2. 
\end{multline*}

Since $\sigma$ is homogeneous function of order $0$ for each variable, $\widetilde \sigma_{j, k}$ given in \eqref{thm-Multiplier-sigma-jk} is independent of $j$ and  $k$, which we now denote by $\widetilde \sigma$. 
We have 
$$
\widetilde \sigma (\xi_1, \xi_2) = \sum_{l_1, l_2 \ge 1} \sigma_{l_1, l_2} \varphi_{l_1}(\xi_1) \varphi_{l_2} (\xi_2) \quad \mbox{ for } \xi_1, \xi_2 \in B_4 \setminus B_{1/4},  
$$
where 
$(\varphi_l)_{l \ge 1}$ is an orthonormal basis of $L^2(B_4 \setminus B_{1/4})$ formed by eigenfunctions of the Laplacian in $B_4 \setminus B_{1/4}$ with the zero Dirichlet boundary condition, i.e., 
$$
-\Delta \varphi_l = \lambda_l \varphi_l \mbox{ in } B_4 \setminus B_{1/4} \quad \mbox{ and } \varphi_l = 0 \mbox{ on } \partial (B_4 \setminus B_{1/4}). 
$$
By the definition, we have, for $l_1, l_2 \ge 1$,  
$$
\sigma_{l_1, l_2} = \int_{B_4 \setminus B_{1/4}} \int_{B_4 \setminus B_{1/4}}  \ws (x_1, x_2) \varphi_{l_1} (x_1) \varphi_{l_2} (x_2) \, dx_1 \, d x_2. 
$$
An integration by parts gives
$$
\sigma_{l_1, l_2} = \frac{1}{\lambda_{l_1}^{N}} \frac{1}{\lambda_{l_2}^{N}} \int_{B_4 \setminus B_{1/4}} \int_{B_4 \setminus B_{1/4}}  \nabla_{x_1}^{N} \nabla_{x_2}^{N} \ws (x_1, x_2) \nabla_{x_1}^{N} \varphi_{l_1} (x_1) \nabla_{x_2}^{N} \varphi_{l_2} (x_2) \, dx_1 \, d x_2, 
$$
which yields, since $N = n+ 2d $ and $\| \nabla^N \varphi_{l} \|_{L^2(B_1 \setminus B_4)} \le C \lambda_{l}^{N/2}$,  
\be \label{thm-Multiplier-p1}
\sup_{l_1, l_2 \ge 1}  \lambda_{l_1}^{n/2+d}  \lambda_{l_2}^{n/2+d} |\sigma_{l_1, l_2}| < + \infty. 
\ee
Here and in what follows in this proof, $C$ denotes a positive constant independent of $f_j$, $j$, $k$, $l_1$, and $l_2$. 

We then have 
$$
\sigma(\xi_1, \xi_2) \hat f_1 (\xi_1)  \hat f_2 (\xi_2)  = \sum_{l_1, l_2 \ge 1} \sigma_{l_1, l_2} \sum_{j, k \in \mZ} \varphi_{l_1}(2^{-j} \xi_1) \varphi_{l_2} (2^{-k}\xi_2) \hat f_{1, j} (\xi_1) \hat f_{2, k} (\xi_2) \quad \mbox{ for } (\xi_1, \xi_2) \in (\mR^d)^2. 
$$
Thus 
\be \label{thm-Multiplier-sum}
\sigma(\xi_1, \xi_2) \hat f_1 (\xi_1)  \hat f_2 (\xi_2) =  \sum_{l_1, l_2  \ge 1} \sigma_{l_1, l_2} \hat g_{1, l_1} (\xi_1)  \hat g_{2, l_2} (\xi_2) \quad \mbox{ for } (\xi_1, \xi_2) \in (\mR^d)^2,
\ee
where, for $\xi_1, \xi_2 \in \mR^d$,  
$$
\hat g_{1, l_1} (\xi_1) = \sum_{j \in \mZ } \varphi_{l_1}(2^{-j} \xi_1)  \hat f_{1, j} (\xi_1) \quad \mbox{ and } \quad \hat g_{2, l_2} (\xi_2) =  \sum_{j \in \mZ } \varphi_{l_2}(2^{-k} \xi_2)  \hat f_{2, k} (\xi_2). 
$$

We thus obtain 
\be \label{thm-Multiplier-sum2}
a (\xi_1, \xi_2) \sigma(\xi_1, \xi_2) \hat \varphi_1 (\xi_1)  \hat \varphi_2 (\xi_2) =  \sum_{l_1, l_2  \ge 1} \sigma_{l_1, l_2} a (\xi_1, \xi_2) \hat g_{1, l_1} (\xi_1)  \hat g_{2, l_2} (\xi_2),
\ee

Applying  \cite[Theorem 7.5.5]{Grafakos2}, we have 
\be \label{thm-Multiplier-p3}
\| T(f_1, f_2) \|_{L^r(\mR^d)} \le C \sum_{l_1, l_2 \ge 1} |\sigma_{l_1, l_2}|  \| g_{1, l_1}\|_{L^{p_1}(\mR^d)} \| g_{2, l_2}\|_{L^{p_2}(\mR^d)}. 
\ee
By the Mihlin-H\"ormander multiplier theorem, see \cite[Theorem 7.9.5]{Hormander1},  one has 
$$
\| g_{1, l_1} \|_{L^{p_1}(\mR^d)} \le  C \|\varphi_{l_1} \|_{H^{n} (B_4 \setminus B_1)} \| f_1\|_{L^{p_1}(\mR^d)} \le C \lambda_{l_1}^{n/2} \| f_1\|_{L^{p_1} (\mR^d)},  
$$
which yields, since $ \|\varphi_{l_1} \|_{H^{n} (B_4 \setminus B_1)} \le C \lambda_{l_1}^{n/2}$, 
\be \label{thm-Multiplier-g1}
\| g_{1, l_1} \|_{L^{p_1}(\mR^d)} \le   C \lambda_{l_1}^{n/2} \| f_1\|_{L^{p_1} (\mR^d)}.   
\ee

Similarly, one obtains 
\be \label{thm-Multiplier-g2}
\| g_{2, l_2} \|_{L^{p_2}(\mR^d)} \le C \lambda_{l_2}^{n/2} \| f_2\|_{L^{p_2}(\mR^d)}. 
\ee

Combining \eqref{thm-Multiplier-p3}, \eqref{thm-Multiplier-g1}, and \eqref{thm-Multiplier-g2} yields 
\be \label{thm-Multiplier-p3b}
\| T(f_1, f_2) \|_{L^r(\mR^d)} \le C \left( \sum_{l_1, l_2 \ge 1} |\sigma_{l_1, l_2}| \lambda_{l_1}^{n/2} \lambda_{l_2}^{n/2} \right) \| f_1\|_{L^{p_1}(\mR^d)} \| f_2\|_{L^{p_2}(\mR^d)}. 
\ee
From  \eqref{thm-Multiplier-p1} and \eqref{thm-Multiplier-p3b}, we derive that  
\be \label{thm-Multiplier-p4}
\| T(f_1, f_2) \|_{L^r(\mR^d)} \le C \left( \sum_{l_1, l_2 \ge 1}  \lambda_{l_1}^{-d} \lambda_{l_2}^{-d} \right)  \| f_1\|_{L^{p_1}(\mR^d)} \| f_2\|_{L^{p_2}(\mR^d)} 
\ee
Since, by the Weyl law, 
$$
\lambda_{l} \ge C l^{2/d}, 
$$
it follows from  \eqref{thm-Multiplier-p4} that 
\be
\| T(f_1, f_2) \|_{L^r(\mR^d)} \le C \sum_{l_1, l_2 \ge 1} l_1^{-2} l_2^{-2}  \| f_1\|_{L^{p_1}(\mR^d)} \| f_2\|_{L^{p_2}(\mR^d)} \le C \| f_1\|_{L^{p_1}(\mR^d)} \| f_2\|_{L^{p_2}(\mR^d)}.
\ee
The proof is complete.  \qed

\section{Proof of \Cref{thm-Multiplier-2}}  \label{sect-thm-Multiplier-2}

Using an appropriate partition of unity of $\partial B \subset (\mR^{d})^m$ where $B$ is the unit ball of $(\mR^{d})^m$, it suffices to prove the result with $\chi \sigma$ instead of $\sigma$ where 
$$
\chi \mbox{ is a function defined on $(\mR^d)^m$, which is homogeneous of degree $0$},  
$$
$$
\chi|_{\partial B} \in C^\infty (\partial B) \mbox{ and } \chi = 0 \mbox{ on } \partial B \setminus B_{\eps} (z) 
$$
for some $z \in \partial B$, where $B_{\eps} (z)$ denotes the ball in $(\mR^d)^m$ centered at $z$ with radius $\eps$. This will be assumed from later on in this proof. 

Since $\eps$ is small, there exists a positive constant $C$ depending only on $m$ and $d$ and there exists $1 \le j \le m$ such that 
$$
|\xi_j | \ge C \max \{ | \xi_1|, \dots | \xi_m |\} \mbox{ for all }   (\xi_1, \dots  \xi_m) \mbox{ in the support of $\chi$}, 
$$
for some positive constant $C$ depending only on $m$ and $d$.  For notational ease, we will take $j=1$. Thus, we may write
$$
|\xi_1 | \ge C \max \{ | \xi_1|, \dots | \xi_m|\}. 
$$

By the multilinear alternating property of $\sigma_m$, we have
$$
\sigma_m (\xi_1, \xi_2, \cdots, \xi_m) = \sigma_m (\xi_1 + \cdots + \xi_m, \xi_2,  \cdots, \xi_m) \quad \mbox{ for } \xi_1, \dots, \xi_m \in \mR^d. 
$$
This implies 
\be
\sigma (\xi_1, \cdots \xi_m) = \mathop{\sum_{1\le \ell_j \le d}}_{\text{for }1 \le j \le k} \Big( (\xi_1 + \cdots \xi_m)_{\ell_1} \cdots (\xi_1 + \cdots + \xi_m)_{\ell_k} C_{\ell_1} \dots C_{\ell_k} \Big),  
\ee
where, with $(e_1, \cdots, e_d)$ being the standard basis of $\mR^d$, 
$$
C_{j} (\xi_1, \dots, \xi_m) = \sigma_m (e_j, \xi_2, \cdots, \xi_m) \quad \mbox{ for } 1 \le j \le d,  
$$
and 
$$
(\xi_1 + \cdots + \xi_m)_l = (\xi_1 + \cdots + \xi_m) \cdot e_l \quad \mbox{ for } 1 \le l \le d. 
$$
Define 
$$
\sigma_{\ell_1, \dots,  \ell_k}(\xi_1, \dots,  \xi_m) = \frac{C_{\ell_1} \cdots C_{\ell_k}}{|\xi_2|^k \cdots  |\xi_m|^k}, 
$$
and 
$$
\ta (\xi_1, \dots,  \xi_m) = \chi (\xi_1, \dots,  \xi_m) \frac{(1+ |\xi_2|^2)^{\frac{k}{2m}} \cdots  (1+ |\xi_m |^2)^{\frac{k}{2m}}}{(1 + |\xi_1|^2)^{\frac{k(m-1)}{2m}}}. 
$$
Set 
\begin{multline}
S (f_1, \dots, f_m) (x)= \int_{(\mR^{d})^d}  \chi (\xi_1, \dots, \xi_m)  \sigma(\xi_1, \dots, \xi_m) \times \\[6pt]
\times \hat f_1 (\xi_1) \cdots \hat f_m (\xi_m) e^{i x (\xi_1 + \cdots +  \xi_m)}   \, d \xi_1 \dots d \xi_m \quad \mbox{ for } x \in \mR^d,  
\end{multline}
and denote 
\begin{multline}
S_{\ell_1 \dots \ell_k} (f_1, \dots, f_m) (x) = 
\int_{(\mR^{d})^d} a (\xi_1, \dots, \xi_m)  \sigma_{\ell_1  \dots \ell_k}(\xi_1, \dots, \xi_m) \times \\[6pt]
\times \hat f_1 (\xi_1) \cdots \hat f_m (\xi_m) e^{i x (\xi_1 + \cdots +  \xi_m)}   \, d \xi_1 \dots d \xi_m \quad \mbox{ for } x \in \mR^d.  
\end{multline}
Set 
$$
\hat g_1 (\xi) =  (1+ |\xi|^2)^{\frac{k(m-1)}{m}}  \hat f_1(\xi),  \quad \hat g_j (\xi) =  \frac{|\xi|^k}{(1+ |\xi|^2)^{k/{2m}}}  \hat f_j(\xi) \quad \mbox{ for } 2 \le j \le m, 
$$
and
$$
\overline{\hat \varphi_{\ell_1 \dots \ell_k}}(\xi) = (\xi_1 + \cdots \xi_m)_{\ell_1} \cdots (\xi_1 + \cdots + \xi_m)_{\ell_k} \overline{\hat \varphi} (\xi), 
$$
We then have 
\be \label{thm-Multiplier-2-c1}
\langle S (f_1, \dots, f_m), \varphi \rangle =  \mathop{\sum_{1\le \ell_j \le d}}_{1 \le j \le k} \langle S_{\ell_1 \dots \ell_k} (g_1, \dots, g_m), \varphi_{\ell_1 \dots \ell_k} \rangle.  
\ee

On the one hand, we have, by \Cref{thm-Multiplier}, 
\be \label{thm-Multiplier-2-c2}
\|S_{\ell_1 \dots \ell_k}(g_1, \dots, g_m) \|_{L^r(\mR^d)} \le C \| g_1 \|_{L^{p_1}(\mR^d)} \cdots \| g_{m}\|_{L^p(\mR^d)}. 
\ee
On the other hand by Mikhlin-H\"ormander's multiplier theorem, we have, for $1 \le j \le m$,  
\be \label{thm-Multiplier-2-c3}
\| g_j\|_{L^{p_j} (\mR^d)} \le C \| f_j \|_{\cL^{p_j}_{s} (\mR^d)}. 
\ee

Combining \eqref{thm-Multiplier-2-c2}, and \eqref{thm-Multiplier-2-c3} yields 
\be
\Big| \langle  S_{\ell_1 \dots \ell_k} (g_1, \dots, g_m), \varphi_{\ell_1 \dots \ell_k} \rangle \Big| \le  C \| f_1 \|_{\cL^{p_j}_{s} (\mR^d)} \cdots \| f_m \|_{\cL^{p_j}_{s} (\mR^d)} \| \varphi\|_{W^{k, r_*}(\mR^d)}, 
\ee
and the conclusion follows from \eqref{thm-Multiplier-2-c1}. \qed

\section{Proof of \Cref{lem-det-H}} \label{sect-lem-det-H}

We set 
$$
\nu_1 = \xi - \xi_1, \quad \nu_2 = \xi_1 - \xi_2, \quad \dots, \quad  \nu_d = \xi_{d-1}. 
$$

Let $\tau$ be a permutation of $\{1, \cdots, d\}$. We have 
\begin{multline}
(-1)^{d} \cF (\det \nabla^2 u) (\xi) = \int_{\mR^{d \times (d-1)}}( \det P_\tau ) \times \\[6pt]
\times \hu (\xi - \xi_1) \hu(\xi_1 - \xi_2)  \cdots \hu(\xi_{d-2} - \xi_{d-1}) \hu(\xi_{d-1}) 
d \xi_1 \dots  \, d  \xi_{d-1}, 
\end{multline}
where the $i$-th column of $P_\tau$ is given by the vector 
$$
( \nu_{\tau(i), i}  \nu_{\tau(i), 1}, \dots, \nu_{\tau(i), i}  \nu_{\tau(i), d} \big)\tr \in \mR^d \quad \mbox{ where } \nu_{j} = (\nu_{j, 1}, \dots, \nu_{j, d})\tr \in \mR^d. 
$$
We have 
$$
( \nu_{\tau(i), i}  \nu_{\tau(i), 1}, \dots, \nu_{\tau(i), i}  \nu_{\tau(i), d} \big)\tr  =  \nu_{\tau(i), i} \nu_{\tau(i)}. 
$$
It follows that 
$$
\det P_\tau = \left( \prod_{1 \le i \le d } \nu_{\tau(i), i} \right) \det \left(\nu_{\tau(1)}, \dots, \nu_{\tau(d)}  \right) = \text{sign} (\tau)\left( \prod_{1 \le i \le d } \nu_{\tau(i), i} \right) \det \left(\nu_1, \dots, \nu_d  \right). 
$$
We derive that 
\begin{multline}
d!(-1)^{d} \cF (\det \nabla^2 u) (\xi) = \int_{\mR^{d \times (d-1)}} \sum_{\tau \text{: permutation}} \text{sign} (\tau)\left( \prod_{1 \le i \le d } \nu_{\tau(i), i} \right) \det \left(\nu_1, \dots, \nu_d  \right) \times \\[6pt]
\times \hu (\xi - \xi_1) \hu(\xi_1 - \xi_2)  \cdots \hu(\xi_{d-2} - \xi_{d-1}) \hu(\xi_{d-1}) 
d \xi_1 \dots  \, d  \xi_{d-1}. 
\end{multline}
This implies 
\begin{multline}
d!(-1)^{d} \cF (\det \nabla^2 u) (\xi) = \int_{\mR^{d \times (d-1)}}  \left(\det \left(\nu_1, \dots, \nu_d  \right) \right)^2 \\[6pt]
\hu (\xi - \xi_1) \hu(\xi_1 - \xi_2)  \cdots \hu(\xi_{d-2} - \xi_{d-1}) \hu(\xi_{d-1}) 
d \xi_1 \dots  \, d  \xi_{d-1}. 
\end{multline}
The conclusion follows from Theorem \ref{thm-Multiplier-2}. \qed 

\begin{remark}[Determinant in 2D]  \label{rem-Hessian} 
\rm It is well-known that in the case $d=2$, the Hessian determinant is continous on $W^{1,2}$. Indeed, one has 
\be \label{d=2}
2 \det (\nabla^2 u) =2 \partial_{12}^2 \Big(\partial_1 u \partial_2 u \Big) - \partial_{11}^2 ( |\partial_{2} u|^2) - \partial_{22}^2 ( |\partial_{1} u|^2). 
\ee
This idendity follows from the usual formula for determiant of gradients 
$$
 \det (\nabla^2 u) = \partial_1 (\partial_1 u C_{11}) +  \partial_2 (\partial_1 u C_{12}) = \partial_1 (\partial_2 u C_{21}) + \partial_2 (\partial_2 u C_{22}),  
$$
where $(C_{ij})$ is the cofactor matrix of $\nabla^2 u$. Thus 
\be \label{d=2-p1}
2 \det (\nabla^2 u) = \partial_1 (\partial_1 u \partial_{22}^2 u) -  \partial_2 (\partial_1 u \partial_{12}^2 u ) -  \partial_1 (\partial_2 u \partial_{21}^2 u) + \partial_2 (\partial_2 u \partial_{11}^2 u),  
\ee
Next we observe that
$$
 \partial_1 (\partial_1 u \partial^2_{22} u)  + \partial_1 (\partial_2 u \partial^2_{21} u)  = \partial^2_{12} \Big( \partial_1 u \partial_2 u \Big), 
$$
$$
\partial_2 (\partial_2 u \partial_{11}^2 u) + \partial_2 (\partial_1 u \partial_{12}^2 u) = \partial_{12}^2 \Big(\partial_1 u \partial_2 u \Big).
$$
and
$$
\partial_2 (\partial_1 u \partial_{12}^2 u) =  \frac{1}{2} \partial_{22}^2 (|\partial_1 u|^2),  \qquad  \partial_1 (\partial_2 u \partial^2_{21} u)  = \frac{1}{2} \partial_{11}^2 ( |\partial_2 u|^2), 
$$
Combining these with   \eqref{d=2-p1}. we obtain the idendity \eqref{d=2}.

\end{remark}

\begin{remark}[Determinant in multiD]  \label{Hessian_formula} 
\rm Iterating the expression \eqref{piola1} for the Jacobian, one can establish that
\[
d (d-1)\; \det (\nabla^2 u) = \sum_{i, j=1}^d \frac{\partial^2}{\partial x_i \partial x_j} \big[\sum_{k\neq i, l \neq j} \frac{\partial u}{\partial x_k} \,  \frac{\partial u}{\partial x_l } \, C_{ij}^{kl} \big]
\]
where the second cofactors \(C_{ij}^{kl}\) of the matrix \( \nabla^2 u \) enjoy symmetry properties, which allow to establish the formula, namely
\[
C_{ij}^{kl} = - C_{kj}^{il} =- C_{il}^{kj}, \qquad  C_{ij}^{kl} = C_{ji}^{lk} .
\]
It is also remarkable that, in dimension larger than 4, the \(C_{ij}^{kl}\) can themselves be expressed again as divergence of lower order multilinear terms.

A more explicit formula is established in \cite{BJ15} 

\[
d! \; \det (\nabla^2 u) = \sum_{\sigma, \tau \in \CS_d} \sgn(\sigma) \sgn( \tau) \frac{\partial^2}{\partial{\sigma(2)}  \partial {\tau(2)} } 
\Big[ \frac{\partial u}{\partial {\sigma(1)}} \, \frac{\partial u}{\partial {\tau(1)}}\,
\frac{\partial^2 u}{\partial{\sigma(3)}  \partial {\tau(3)} }\, \cdots \, \frac{\partial^2 u}{\partial{\sigma(d)}  \partial {\tau(d)} } \Big].
\]

\end{remark}

\bigskip 
\noindent {\bf Acknowledgement:}  The authors wish to thank Nicolas Lerner for many fruitful discussions.


\end{document}